\documentstyle{amsppt}
\magnification=\magstep1
\vsize=21.1 true cm
\hsize=15.5 true cm
\voffset=0 true cm
\pageno=1
\nologo
\topmatter
\title{Non-zero degree maps between $2n$-manifolds}
\endtitle

\author{Haibao Duan and Shicheng Wang}
\endauthor
\affil Institute of Mathematics, \\
Chinese Academy of Science,
Beijing 100080, China\\
dhb$\@$math03.math.ac.cn\\
and\\
 LMAM, Department of Mathematics, \\
 Peking University, Beijing 100871, China\\
 wangsc$\@$math.pku.edu.cn
\endaffil
\thanks{{Both author are supported by MSTC, NSFC.
The comments of F. Ding, J.Z. Pan, Y. Su and the referee enhance
the paper.}}
\endthanks
\email
\endemail
\abstract{Thom-Pontrjagin constructions are used to give a
computable necessary and sufficient condition when a homomorphism
$\phi : H^n(L;Z)\to H^n(M;Z)$ can be realized by a map $f:M\to L$
of degree $k$ for closed $(n-1)$-connected $2n$-manifolds $M$ and
$L$, $n>1$. A corollary is that each $(n-1)$-connected
$2n$-manifold admits selfmaps of degree larger than 1, $n>1$.

In the most interesting case of dimension 4, with the additional
surgery arguments we give a necessary and sufficient condition for
the existence of a degree $k$ map from a closed orientable
4-manifold $M$ to a closed simply connected 4-manifold $L$ in
terms of their intersection forms, in particular there is a map
$f:M\to L$ of degree 1 if and only if the intersection form of $L$
is isomorphic to a direct summand of that of $M$.

\vskip 0.5 true cm

\vskip 0.5 true cm


}

\endabstract
\endtopmatter
\document
%\normalbaselineskip=1.44\normalbaselineskip
%\normalbaselines

{\bf \S 1. Introduction and Results.}

A fundamental question in manifold topology is the following

\proclaim{Question} For given closed oriented manifolds $M$ and
$L$ of the same dimension, can one decide if there is a map
$f:M\to L$ of degree $k\ne 0$? in particular $k=1$?
\endproclaim

For dimension 2 the answer is affirmative and simple (see [1]).
But in general this question is difficult. For dimension 3 many
rather general facts have been discussed in last  ten years (see a
survey paper [2]). For dimension at least 4, there are only some
special examples discussed up to our knowledge (see references in
[3]). In this paper we will obtain affirmative answers for maps
between $(n-1)$-connected $2n$-manifolds, $n>1$, and for maps from
arbitrary 4-manifolds to simply connected 4-manifolds.

Unless special indications are made, in this paper all
$2n$-manifolds are closed, connected and orientable and  $n>1$.
All matrices, homology rings and cohomology rings are in integers.

Let $M$ be a $2n$-manifold
and let $\bar H^n(M)$ (resp. $\bar H_n(M)$)
be the free part of $H^n(M)$ (resp. $H_n(M)$).
Then the cup product operator
$$ \bar H^n(M)\otimes \bar H^n(M)\to H^{2n}(M)\tag 1.1$$
defines the intersection form $X_M$ over $\bar H^n(M)$,
which is bi-linear and unimodular.

Let $\alpha^*=(a^*_{1},...,a^*_{m})$ be a basis for $\bar H^n(M)$.
Then $X_M$ determines an $m$ by $m$ matrix $A=(a_{ij})$, where
$a_{ij}$ is given by $a^*_{i}\cup a^*_{j}=a_{ij}[M]$, and $[M]$ is
the fundamental class of $H^{2n}(M)$. We often present this fact
by the formula
$${\alpha^*}^\bot\cup\alpha^*=A[M], \tag 1.2$$
where ${\alpha^*}^\bot$ is the transpose matrix of $\alpha^*$.
We also have

$${\alpha_*}^\bot\cap\alpha_*=A, \tag 1.3$$
where $\alpha_*\subset \bar H_n(M)$ is the Poincar\'e dual basis
of $\alpha^*$. (We use $\cap$ to denote the intersection of
homology classes in this paper.)

Suppose  $M$ is a closed $(n-1)$-connected $2n$ manifold. Then $M$
has a presentation $V_{i=1}^m S^n_i \cup_g D^{2n}$, i.e., $M$ is
obtained from $V_{i=1}^m S_i^n$, the $m$-fold wedge of $n$-spheres
with one point $y$ in common, by attaching a $2n$-cell $D^{2n}$
via an attaching map
$$g: S^{2n-1}=\partial D^{2n}\to V_{i=1}^m S_i^n,$$
and the homotopy type of $M$ is determined by $g\in
\pi_{2n-1}(V_{i=1}^m S_i^n)$. Such a presentation of $M$ provides
$M$ a complete set of homotopy invariants $(A ;t_1,...,t_m) $,
where  $t_i\in \pi_{2n-1}(S^n)$, and $A_{m\times m}$ is the
linking matrix. (For details, including the relation between the
linking matrix and the intersection matrix, see Section 2).

Suppose $L$ is another closed $(n-1)$-connected $2n$-manifold with
presentation $V_{j=1}^l S_j^n\cup_h D^{2n}$ and $h$ corresponding
to $(B;u_1,...,u_l)$.

For each map $f:M\to L$, we use $f_*$ and $f^*$ to denote the
induced homomorphisms on the homology rings and cohomology rings.
Choose bases $\alpha=\{S^n_i, i=1,...,m\}$ and $\beta=\{S^n_j,
j=1,...,l\}$  of $H_n(M)$ and $H_n(L)$ respectively. Let
$P=(p_{ij})$ be an $m$ by $l$ matrix.

Now we can state the following Theorem 1 which gives a necessary
and sufficient condition for the existence of a degree $k$ map
between $(n-1)$-connected 2n-manifolds $M$ and $L$, $n>1$.

\proclaim{Theorem 1} Suppose $M$ and $L$ are  closed oriented
$(n-1)$-connected $2n$-manifolds with presentations and data as
above. Then there is a map $f: M\to L$ of degree $k$ such that
$f_* (\alpha) =\beta P^\bot$ if and only if $P$  satisfies the
following equations

$$P^\bot A P=kB \tag 1.4$$
and
$$ku_r=\sum_v p_{vr}t_v +(\sum _v \frac 12 p_{vr}(p_{vr}-1)a_{vv}+
\sum _{v<w}  p_{vr}p_{wr}a_{vw})[s^n,s^n], \tag 1.5$$

where $s^n$ generates  $\pi_n(S^n)$ and
$[s^n, s^n]\in \pi_{2n-1}(S^n)$ is the Whitehead product.

Moreover the homotopy classes of maps $f: M\to L$ such that $f_*
(\alpha) =\beta P^\bot$ are 1-1 corresponding to the elements of
(the finite group) $\pi_{2n}(L)$.
\endproclaim

The following Theorem 2 is stated for all maps between two closed
oriented $2n$-manifolds, and its proof is direct.

\proclaim{Theorem 2} Suppose $M$ and $L$ are  closed oriented
$2n$-manifolds with  intersection matrices $A$ and $B$ under some
given bases $\alpha^*$ for $\bar H^n(M)$ and $\beta^*$ for $\bar
H^n(L)$ respectively. If there is a map $f: M\to L$ of degree $k$
such that $f^* (\beta^*) =\alpha^* P$, then
$$P^\bot AP=kB.\tag 1.6$$
Moreover if $k=1$, then $X_L$ is isomorphic to a direct summand of $X_M$.
\endproclaim

In dimension 4, we have the following necessary and sufficient
condition for the existence of a degree $k$ map from a closed
oriented 4-manifold to a simply connected 4-manifold.

\proclaim{Theorem 3} Suppose $M$ and $L$ are closed oriented
4-manifolds with  intersection matrices  $A$ and $B$ under  given
bases $\alpha^*$ for $\bar H^2(M)$ and $\beta^*$ for $\bar H^2(L)$
respectively. If $L$ is simply connected, then there is a map $f:
M\to L$ of degree $k$ such that $f^* (\beta^*) =\alpha^* P$ if and
only if
$$P^\bot A P=kB\tag 1.7$$

Moreover there is a map $f:M\to L$ of degree 1
if and only if
$X_L$ is isomorphic to a direct summand of $X_M$.
\endproclaim

Note that $-X_M=X_{-M}$. So similar results can be stated for maps
of degree $-1$.

\proclaim{Corollary 1}  Suppose  $f:M\to L$ is a map of degree 1
between two simply connected oriented 4-manifolds. Then $M$ has
the  homotopy type  of $L\# Q$, where $Q$ is a simply connected
4-manifold; Moreover if $\chi(M)=\chi(L)$, then $M$ and  $L$ are
homotopy equivalent (or homeomorphic if $X_M$ is even).
\endproclaim

\vskip 0.5 true cm

{\bf Definition:}
For any two closed oriented $n$-manifolds $M$ and $L$, define
$$D(M,L)= \{deg f | f: M\to L\},\tag 1.8$$
that is, $D(M,L)$ is the set of all degrees of maps from $M$ to $L$.

\proclaim{Corollary 2} Suppose $M$ is an $(n-1)$-connected
$2n$-manifold, $n>1$. Let $T(n)$ be the order of the torsion part
of $\pi_{2n-1}(S^n)$. Then $k^2\in D(M,M)$ if $k$ is a multiple of
$2T(n)$ when $T(n)$ is even,  or a multiple of $T(n)$ otherwise.
Moreover if $M$ is a simply connected 4-manifold, then $D(M,M)$
contains all squares of integers.
\endproclaim

\proclaim{Corollary 3} Suppose $M$ is a closed oriented
$2n$-manifold and $\text{dim}\bar H^n(M)=m$ is odd. Then for any
map $f: M\to M$, the degree $k$ of $f$ is a  square of an integer.
\endproclaim

Say the intersection form $X_M$ is even if for any $x\in \bar H_n(M)$,
$x\cap x$ is even.

\proclaim{Corollary 4} Let $f:M\to L$ be a map of degree $k$
between closed oriented $2n$-manifolds. If $X_M$ is even
and $X_L$ is not, then $k$ must be even.
\endproclaim

{\bf Definition:} Let $M, L$ be two $n$-dimensional closed
oriented manifolds.
 Say $M$ dominates $L$
if there is a map $f:M\to L$ of non-zero degree.

\proclaim {Corollary 5} The homotopy types of $(n-1)$-connected
$2n$-manifolds dominated by any given closed oriented
$2n$-manifold $M$ are finite. Each closed oriented 4-manifold
dominates at most finitely many simply connected 4-manifolds.
\endproclaim

The paper is organized as follow: In section 2 we discuss related
facts about Thom-Pontrjagin constructions and the presentations of
$(n-1)$-connected $2n$-manifolds, which will be used
 to prove Theorem 1 in Section 3.
Theorem 2 and Theorem 3 are proved in Section 4 and Section 5
respectively, the corollaries are proved in Section 6. In Section
7 we give some applications to concrete examples in dimension 4.

We  end the introduction by comparing this paper with [3], another
paper of the same authors. The major part of Theorem 1 in this
paper is the same as Theorem A in [3], which is proved in this
paper by using Thom-Pontrjagin construction, a rather elementary
and geometric argument, as presented below, and is proved in [3]
by more advanced algebraic topology argument. The remaining parts
of the two papers are different: The present paper keeps the
geometric style and is devoted mainly on maps from 4-manifolds to
simply connected 4-manifolds, and however [3] carries the
advantage of modern algebraic topology and is devoted to the
computation of $D(M,N)$ for some classes of $(n-1)$-connected
$2n$-manifolds, $n>2$.

\vskip 0.5 true cm

{\bf \S 2. Thom-Pontrjagin construction and presentation of
$(n-1)$-connected $2n$-manifold.}

All results in this section must be known.  But we can not find
the statements of Theorem 2.2 and Lemma 2.4  below in the
literature, and therefore we provide proofs for them.

Let $f: S^{2n-1}\to S^n$ be any map with $y\in S^n$ a regular
value. Then $W=f^{-1}(y)\subset S^{2n-1}$, an (n-1)-manifold,  and
$\tau =f^*(T_y S^n)$, the pull-back of the tangent space of $S^n$
at $y$, together provide a framed normal bundle of $W$ in
$S^{2n-1}$. Call the pair $(W, \tau)$ a manifold in $S^{2n-1}$
with a given framed normal bundle.

If both $S^{2n-1}$ and $S^n$ are oriented, then both $W$ and the
normal framed bundle are naturally oriented as follows: The normal
space at each point of $W$ has the orientation lifted from the
orientation of the tangent space at $y\in S^n$, and $W$ is
oriented so that the orientations of $W$ and its normal bundle
give the orientation of $S^{2n-1}$.

On the other hand, for any framed normal bundle $(W, \tau)$, the
$\epsilon$-neighborhood $W_\epsilon$  of $W$ in $S^{2n-1}$ has a
product structure $D^n\times W $ induced by the given framing. Let
$p: D^n\times W\to D^n$ be the projection, $q: D^n\to D^n/\partial
D^n=S^n$ be the quotient map, then we can extend $q\circ p:
W_\epsilon \to S^n$ to a map $f: S^{2n-1}\to S^n$ by sending the
remaining part of $S^{2n-1}$ to the point $q(\partial D^n)\in
S^n$.

Define

$$\Cal F(n)=\{\text{ $(n-1)$-manifold $W^{n-1}\subset S^{2n-1}\ $
with a framed normal bundle}\}$$

Note the underlying space of an element in $\Cal F(n)$ does not
need to be connected.

Let $\Omega(n)$ be the framed bordism classes of $\Cal F(n)$.
For each $s\in \Cal F(n)$,
we use $[s]$ for its framed bordism class.

\proclaim{Theorem 2.1} The construction (Thom-Pontrjagin) above
gives a 1-1 and onto map $T:  \Omega(n) \to \pi_{2n-1}(S^n)$. (see
[p.209 4] for a proof).
\endproclaim

For each $s=(W, \tau)\in \Cal F(n)$, we will use
$-s$ for $(-W, \tau)$.
For $s'=(W',\tau')\in \Cal F(n)$, which is disjoint from $s$,
let $lk(s, s')=lk(W, W')$,
and $s\sqcup s'$ be the disjoint union of $s$ and $s'$, which is
again in $\Cal F(n)$.

\proclaim{Theorem 2.2}
For any two disjoint $s_1, s_2\in\Cal F(n)$,
$$T([s_1\sqcup s_2])=
T([s_1])+T([s_2])+lk (s_1, s_2)[s^n, s^n],\tag 2.1$$
where $s^n$  generates $\pi_n(S^n)$ and
$[s^n, s^n]\in \pi_{2n-1}(S^n)$ is the Whitehead product.
\endproclaim

A result similar to Theorem 2.2 is proved  in [pp167-168 Lemma 2,
5], for diffeotopic classes (which are  finer than framed bordism
classes) of $n-$spheres (which form a proper subset of
$n-$manifolds). Indeed Theorem 2.2 can be derived from [Lemma 2,
5] by first a projection and then some modifications. Since the
statement of [Lemma 2, 5] involves many terminologies in homotopy
theory and its proof also involves deep results of Smale and
Kervaire, we prefer to give a direct proof of Theorem 2.2 as
below.

\demo{\bf Proof}
We first verify that

$$T([s_1\sqcup s_2])=T([s_1])+T([s_2])\qquad \text{if} \qquad lk(s_1,s_2)=0.
\tag 2.2$$

If $s_1$ and $s_2$ are in the upper-hemisphere $S_+^{2n-1}$ and
lower-hemisphere $S_-^{2n-1}$ respectively, then clearly
$T([s_1\sqcup s_2])=T([s_1])+T([s_2])$. In general, by the
condition $lk(s_1,s_2)=0$ and standard arguments in geometric
topology, one can construct a framed bordism $W_i \subset
S^{2n-1}\times [0,1]$ between $s_i \subset S^{2n-1}\times \{0\}$
and $s_i^*\subset S^{2n-1}\times \{1\}$, $i=1,2$, such that

(1) $W_1$ and $W_2$ are disjoint in $S^{2n-1}\times [0,1]$,

(2) $s_1^*$ and $s_2^*$ are in  $S_+^{2n-1}\times \{1\}$ and
$S_-^{2n-1}\times \{1\}$ respectively.

Then it follows that (2.2) still holds.

Now suppose $lk(s_1, s_2)=k\ne 0$.

First pick a $(2n-1)$-ball $B^{2n-1}$ in $S^{2n-1}$ which is disjoint from
$s_1\sqcup s_2$.
Then pick two oriented $n$-balls
$B_1, B_2\subset \text{int}B^{2n-1}$,
with $lk(\partial B_1, \partial B_2)=k$.

Let  $s_i'=(\partial B_i, \tau_i') \in \Cal F(n)$,
where the normal framing
$\tau_i'$ can be extended to $B_i$, $i=1,2$.
Clearly we have
$$[s_1']= [s_2']=0\in \Omega (n), \tag 2.3$$

and

$$lk(\pm s'_1, s'_2)=\pm k,\qquad lk(s_i, s_j')=0 \qquad
 i,j=1,2. \tag 2.4$$

Then by (2.2), (2.3) and (2.4), we have

$$T([s_1\sqcup s_2])+T([-s_1'\sqcup s_2'])$$
$$=T([(s_1\sqcup s_2)\sqcup (-s_1'\sqcup s_2')])=
T([(s_1\sqcup -s_1')\sqcup (s_2\sqcup s_2')])$$
$$=T([s_1\sqcup -s_1'])+T([s_2\sqcup s_2'])=T([s_1])+T([s_2]).\tag 2.5$$

From the facts that $[s_1']=[s_2']=0\in \Omega(n)$, and their base
spaces
 are the boundaries of
$n$-balls,  and the definitions of Thom-Pontrjagin construction
and of Whitehead product [pp138-139, 6], one can verify  that
$T([s_1'\sqcup s_2'])=lk (s_1', s_2')[s^n, s^n]$. Hence
$T([-s_1'\sqcup s_2'])=-lk (s_1, s_2)[s^n, s^n]$. Substitute it
into (2.5). We finish the proof of Theorem 2.2. \qed

\vskip 0.3 true cm

Call $M$ a pre-n-space, if $M$ has the homotopy type of a space
obtained from $V_{i=1}^m S_i^n$, the $m$-fold wedge of $n$-spheres
with one point $y$ in common, by attaching a $2n$-cell $D^{2n}$
via an attaching map
$$g: S^{2n-1}=\partial D^{2n}\to V_{i=1}^m S_i^n.\tag 2.6$$
So the homotopy type of $M$ is determined by
$g\in \pi_{2n-1}(V_{i=1}^m S_i^n)$.
For short denote this space by
$V_{i=1}^m S^n_i \cup_g D^{2n}$.

Suppose both $S^n_i$ and $D^{2n}$ are oriented. Let $y_i\in S_i^n$
be a point (other than the base point $y$). Make $g$ transverse to
$y_i$ for each $i=1,...,m$ so that $K_i=g^{-1}(y_i)$ is an
oriented $(n-1)$-manifold with a framed normal bundle in
$S^{2n-1}$. The orientation of $K_i$ is given by the orientations
of $\partial D^{2n}= S^{2n-1}$ and $S_i^n$, and the orientation of
the framed normal bundle of $K_i$ is given by Thom-Pontrjagin
construction, i.e., the pull-back of the orientation of the
tangent space at $y_i$ of $S^n_i$. The oriented framed
(n-1)-manifolds $K=\{K_1,...,K_m\}$ give a linking matrix
$A=(a_{vw})$, where $a_{vw}=lk(K_v, K_w)$, the linking number of
$K_v$ and $K_w$. The manifold $K_i$ with the given framed normal
bundle gives an element $t_i\in \pi_{2n-1}(S^n)$.

So we have a set of invariants $(A; t_1,...,t_m)$ of $M$. Denote
$$\Cal I(m,n)=\{(A_{m\times m};t_1,...,t_m \in \pi_{2n-1}(S^n))\},\tag 2.7$$
where $A$ is symmetric if $n$ is even and is anti-symmetric if $n$
is odd. Moreover $M$ here has  the homotopy type of a manifold,
thus $A$ is unimodular.

\proclaim{Theorem 2.3} There is a 1-1 correspondence $H:
\pi_{2n-1}(V_{i=1}^m S_i^n)\to \Cal I(m,n)$. Hence

(1) $(A;t_1,...,t_m)$
is a complete homotopy type invariant
of the pre-n-space $M$.

(2) Moreover if $n=2$, then $A$ itself is a complete
homotopy type invariant.
\endproclaim

\demo{\bf Proof}  The result is known. For (1) see [p. 128, 5].
For (2) see [p. 23, 7], and (2) is also derived from (1) by the
fact that $\pi_3(S^2)=Z$ is torsion free and hence $a_{ii}$
determines $t_i$. The proof is an  application of Thom-Pontrjagin
construction. See [p. 23, 7].\qed

\proclaim{Lemma 2.4} The linking matrix $A$ provided by presentation
$V_{i=1}^m S^n_i \cup_g D^{2n}$ is  the intersection matrix
under the algebraic dual basis $\alpha^*$ of
$\alpha=([S_{1}^n],...,[S_{m}^n])$.
\endproclaim

\demo{\bf Proof} Now we can further homotope $g$ above so that
$K_i=g^{-1}(y_i)\subset S^{2n-1}$ is connected. Let $F_i'\subset
D^{2n}$ be an oriented n-chain bounded by $K_i$. Let $F_i$  be the
cycle in $M=V_{i=1}^m S^n_i \cup_g D^{2n}$, defined by $F_i'/K_i$,
where $g(K_i)=y_i$. Then

(a) $S_i^n\cap F_j$ is a point if $i=j$ and is empty if $i\ne j$,

(b) the intersection matrix of $M$ defined by cap product under
the basis $\alpha_*=(F_1,...,$ $F_m)\subset \bar H_n(M)$ is the
linking matrix $A$ provided by $g$.

Let $\alpha^*$ be the Poincar\'e dual basis of $\alpha_*$. (b)
implies that the intersection matrix of $M$ under the basis
$\alpha^*$ is $A$ and then  (a) implies that $\alpha^*$ is the
algebraic dual of $\alpha$.\qed

\vskip 0.5 true cm
{\bf \S 3. Proof of Theorem 1.}

Recall that $M$ and $L$ are  closed $(n-1)$-connected
$2n$-manifolds with presentations

$$M=V_{i=1}^m S^n_i \cup_g D^{2n}\qquad \text{and}
\qquad L=V_{j=1}^l S_j^n\cup_h D^{2n}.\tag 3.1$$
We can consider
$$\alpha=\{S^n_i, i=1,...,m\}, \qquad
\text{and}\qquad  \beta=\{S^n_j, j=1,...,l\}\tag 3.2$$ as bases
for either $H_n(M)$ and $H_n(L)$, or $H_n(V_{i=1}^m S^n_i)$ and
$H_n(V_{j=1}^l S_j^n)$.

\proclaim{Proposition 3.1} Let $f : V_{i=1}^m S_i^n \to V_{j=1}^l
S_j^n$ be a map such that $f_* \alpha  =\beta P^\bot$. Then for
the element $g\in \pi_{2n-1}(V_{i=1}^m S^n_i)$ corresponding to
$(A;t_1,...,t_m)\in \Cal I(m,n)$, the element $f_*(g)\in
\pi_{2n-1}(V_{j=1}^l S^n_j)$ is corresponding to
$(A';t'_1,...,t'_l)\in \Cal I(l,n)$, where
$$A'=P^\bot A P \tag 3.3$$
and

$$t'_r=\sum_v p_{vr}t_v +(\sum _v \frac 12 p_{vr}(p_{vr}-1)a_{vv}+
\sum _{v<w}  p_{vr}p_{wr}a_{vw})[s^n,s^n]. \tag 3.4$$
\endproclaim

\demo{\bf Proof} Call a map from an $n$-disc $D^n$ to an
$n$-sphere $S^n$ with base point $z$ a boundary pinch if the
interior of $D^n$ is mapped onto $S^n-z$ homeomorphically and
$\partial D^n$ is mapped onto $z$. If both $D^n$ and $S^n$ are
oriented, we can say a boundary pinch is of degree 1 or degree
$-1$.

Recall  $P=(p_{ij})_{m\times l}$, then
$P^\bot =(p_{ji})_{l\times m}$.

For each $i$, take $\sum_j |p_{ij}|$ disjoint $n$-discs in $S_i^n$
away from the common wedge point $y \in V_{i=1}^m S_i^n$. Each
$n$-disc has the induced orientation of $S^n_i$. Then map
$|p_{ij}|$ discs to $S_j^n$ via boundary pinches of degree 1 if
$p_{ij}$ is positive or $-1$ otherwise. Doing this for all $i$,
and finally mapping the complement of all these $\sum _{i,j}
|p_{ij}|$ discs to the base point $z\in V_{j=1}^l S_j^n$, we get a
map $f': V_{i=1}^m S_i^n \to V_{j=1}^l S_j^n$  and clearly
$f^{\prime}_* \alpha =\beta P^\bot$ on $H_n=\pi_n$. Since the
homotopy group $\pi_i$ vanishes for $0<i<n$, up to homotopy, we
may assume that $f=f'$.

Let $z_j\in S_j^n-z$. Note $f^{-1}(z_j)$ consists of $\sum_i
|p_{ij}|$ points, and $|p_{ij}|$ points of them lie in $S_i^n$
with sign 1 if $p_{ij}$ is positive or sign $-1$ otherwise.
Suppose $g$ is transverse to $y_i\in S^n_i -y$. Then the oriented
$(n-1)$-manifold $K_i=g^{-1}(y_i)$ and its lifted framed normal
bundle $g^*(T_{y_i}S_i^n)$ provide an element $s_i\in \Cal F(n)$.
Moreover $A=(lk (K_i, K_j))$ and $t_i=T([s_i])$ provide the set of
invariants $(A;t_1,...,t_m)$. Now suppose $f\circ g$ is transverse
to $z_r$ and

$$(f\circ g)^{-1}(z_r)=g^{-1}(f^{-1}(z_r))=\sqcup _{v=1}^m
\sqcup _{e=1}^{|p_{vr}|} \text{sign($p_{vr}$)} K_v^e, \tag 3.5$$
where $K_v^e$ is the preimage of a point in $f^{-1}(z_r)\cap
S^n_v$ under $g$. Let $s_v^e\in \Cal F(n)$ be $K_v^e$ with its
lifted framed normal bundle. We have

$$t'_r=T([\sqcup _{v=1}^m \sqcup_ {e=1}^{|p_{vr}|}
\text{sign($p_{vr}$)} s_v^e]).\tag 3.6$$

Then
$$lk (s_v^e, s_w^{e'})=
lk (K_v^e, K_w^{e'})=lk (K_v, K_w)=a_{vw}\,\, \text {and}\,\,
T ([s_v^e])=T([s_v])=t_v.\tag 3.7$$

Now we use abbreviations as below
$$
p_{vr} K_v= \sqcup_ {e=1}^{|p_{vr}|} \text{sign}(p_{vr}) K_v^e,
\qquad p_{vr} s_v= \sqcup _{e=1}^{|p_{vr}|} \text{sign}(p_{vr})
s_v^e.\tag 3.8$$

Now (3.3) follows from (3.5), (3.7), (3.8) and definition, since
$$A'=(a^{\prime}_{rs})=(lk ((f \circ g)^{-1}(z_r),(f\circ g)^{-1}(z_s)))
=(lk( \sqcup _v p_{vr} K_v, \sqcup _w p_{ws} K_w))$$
$$= (\sum _{v, w} p_{vr}p_{ws}lk (K_v, K_w))
= (\sum _{v, w} p_{vr}p_{ws} a_{vw})=P^\bot A P. \tag 3.9$$

By (3.6), (3.7), (3.8), and by applying Theorem 2.2 to $T([\sqcup
_v p_{vr} s_v])$ and $T([p_{vr}s_v])$, (3.4) follows, since

$$T([\sqcup _v p_{vr} s_v])
=\sum _v T([p_{vr}s_v])+\sum_{v<w} p_{vr}p_{wr}lk (s_v,
s_w)[s^n,s^n]$$
$$=\sum_v p_{vr}T([s_v]) +(\sum _v \frac 12 p_{vr}(p_{vr}-1)a_{vv}+
\sum _{v<w}  p_{vr}p_{wr}a_{vw})[s^n,s^n]. \tag 3.11$$

Proposition 3.1 is proved.
\qed
\demo{\bf Proof of Theorem 1}

{\bf (1) The sufficient part.} We are going to show that if the
complete homotopy invariants of $M$ and $L$,  and the matrix $P$
satisfy equations (1.4) and (1.5), then there is a map $f:M\to L$
of degree $k$ such that $f_* \alpha=\beta P^\bot$.

Let $J=(J_1,...,J_l)=h^{-1}(z_1,...,z_l)\subset B^{2n-1}\subset
S^{2n-1}$ be the (n-1)-manifold provided by $h:S^{2n-1}\to
V_{j=1}^l S_j^n$ and $B^{2n-1}$ is a (2n-1)-ball. Let  $p_k:
S^{2n-1}\to S^{2n-1}$ be a branched covering of degree $k$ over a
standard sphere $S^{2n-3}\subset S^{2n-1}$ which misses
$B^{2n-1}$. Then $p_k^{-1}(J)=(J^1,...,J^k)$, each
$J^i=(J^i_1,...,J^i_l)$ lies in a homeomorphic lift $B^{2n-1}_i$
of $B^{2n-1}$. It follows that

$$lk(J^i_v, J^i_w)=lk (J_v, J_w)\,\, \text{and}\,\, lk(J^i_v,J^j_w)=0\,
\text{if $i\ne j$}.\tag 3.12$$

Now consider the composition
$$
\CD S^{2n-1}@> p_k>> S^{2n-1} @> h>> V_j S_j^n.
\endCD                              \tag 3.13
$$
By (3.12), the discussion above and Theorem 2.2, it is easy to
verify that the element in $\Cal I(l,n)$ corresponding to $[h\circ
p_k]=k[h]\in \pi_{2n-1}(V_{j=1}^l S_j^n)$ is $(kB;
ku_1,...,ku_l)$.

Let $f: V_{i=1}^m S_i^n\to V_{j=1}^l S^n_j$ be a map such that
$f_*(\alpha)=\beta P^\bot$. Then by (3.3) and (3.4) of Proposition
3.1, we have
$$[f \circ g] =[h\circ p_k], \tag 3.14$$
that is

$$f_*[g] =k[h]\in \pi_{2n-1}(V_{j=1}^l S^n_j). \tag 3.15$$
By (3.15), $f: V_i S_i^n \to  V_j S_j^n$ can  be extended to a
map, still denoted by $f$,

$$M=V_{i=1}^m S^n_i \cup_g D^{2n} \to N=V_{j=1}^l S^n_j \cup_h D^{2n}.
\tag 3.16$$

Clearly $f_*(\alpha)=\beta P^\bot$. It is also easy to verify that $f$ is of
degree $k$, since $p_k$ is of degree $k$.
\qed

{\bf (2) The necessary part.} Suppose there is a map $f:M\to L$ of
degree $k$ such that $f_* \alpha=\beta P^\bot$.  We are going to
show that the complete homotopy invariants of $M$ and $L$ and the
matrix $P$ satisfy equations (1.4) and (1.5).

Up to homotopy, we may assume that $f(V_i S_i^n)\subset V_j S_j ^n
$ and the restriction on $V_i S_i^n$ is the same as the map $f'$
given in the proof of Proposition 3.1. Then by Proposition 3.1,
the invariant $(A';t_1',...,t_l')$ associated with the map $f\circ
g$ is given by (3.3) and (3.4).

Since the kernel of

$$i_*: \pi_{2n-1}(V_j S_j^n)\to \pi_{2n-1}(V_j S_j^n\cup_h D^{2n})
=\pi_{2n-1}(V_j S_j^n)/[h] \tag 3.18$$

induced by the inclusion $i : V_j S_j^n\to V_j S_j^n\cup_h D^{2n}$
is the cyclic group $<[h]>$, and $[i\circ f\circ g]=0$, we have
$[f\circ g]=d[h]$ for some integer $d$. It follows that (1.4) and
(1.5) hold if we replace $k$ by $d$.

Under the algebraic dual bases $\alpha^*$ and $\beta^*$ of
$\alpha$ and $\beta$, the intersection matrices are $A$ and $B$
for $M$ and $L$ respectively by Lemma 2.4. Moreover $f^*\beta^* =
\alpha^* P$. Since $f$ is of degree $k$, we have $P^\bot A P=kB$
by Theorem 2 (note the proofs of Theorem 2 and Lemma 2.4 are
independent of Theorem 1). So we have $d=k$. \qed

{\bf (3) The 1-1 correspondence.} Let $\Cal P=\{f:M\to L| f_*
\alpha = \beta P^\bot\}$. Since $M$ and $L$ are $(n-1)$-connected
$2n$-manifolds, for any $f\in \Cal P$, we may assume that up to
homotopy
$$f', f: M=V_{i=1}^m S^n_i \cup_g D^{2n} \to N=V_{j=1}^l S^n_j \cup_h D^{2n},$$
differ only in the interior of $D^{2n}$, where $f'$ is constructed
in the proof of Proposition 3.1. Now the restrictions of $f$ and
$f'$ on $D^{2n}$ give a map $C(f,f'): S^{2n}\to L$, which presents
an element in $\pi_{2n}(L)$.   Then the verification that
$f\mapsto C(f,f')$ for $f\in \Cal P$ induces a 1-1 correspondence
$\{[f]|f\in \Cal P\} \to \pi_{2n}(L)$ is routine.

(The finiteness of $\pi_{2n}(L)$ follows from [p. 45, Kahn, Math.
Ann., Vol. 180]).\qed

\vskip 0.5 true cm
{\bf \S 4. Proof of Theorem 2.}

Recall we have
$${\alpha^*}^\bot\cup \alpha^* =A[M]=(a_{ij})[M],\qquad {\beta^*}^\bot \cup
\beta^*=B[L]=(b_{ij})[L].\tag 4.0$$

Suppose $f: M\to L$ is of degree $k$ such that $f^* (\beta^*)
=\alpha^* P$. Then we have

$$f^*([L])=k[M].\tag 4.1$$

By (4.0)  and (4.1), we have

$$kB[M]=(b_{ij}k[M])=(f^*(b_{ij}[L]))=
(f^*(b_{i}^*\cup b_{j}^*))$$
$$=(f^*(b_{i}^*)\cup f^*(b_{j}^*))=(f^*(\beta^*))^\bot\cup (f^*\beta^*)=
(\alpha^* P)^\bot \cup \alpha^* P\tag 4.2$$
$$=P^\bot ({\alpha^*}^\bot \cup \alpha^*) P=P^\bot A P[M].$$

That is

$$P^\bot A P=kB.\tag 4.3$$

For any $b_i^*, b_j^* \in \beta^* \subset \bar H^n(L)$, we have
$$X_M(f^*(b_i^*), f^*(b_j^*))[M]=f^*(b_i^*)\cup f^*(b_j^*)$$
$$=f^*(b_i^*\cup b_j^*)=f^*(b_{ij}[L])=b_{ij} f^*([L])=b_{ij}[M].\tag 4.4$$
Hence the restriction of $X_M$ on the subgroup $f^*(\bar
H^n(L))\subset \bar H^n(M)$ is isomorphic to $X_L$.  In particular
it is still unimodular. By orthogonal decomposition lemma [p.5,
8],
$$X_M=X_{f^*{\bar H^n(L)}}\oplus X_{H'}=X_L\oplus X_{H'},$$

where $H'$ is the orthogonal complement of $f^*{\bar H^n(L)}$ and
$X_{H'}$ is the restriction of $X_M$ on $H'$.
 \qed

\vskip 0.5 true cm

{\bf \S 5. Proof of Theorem 3.}

\proclaim{Proposition 5.1} Let $M$ and $L$ be two
closed oriented 4-manifolds satisfying
$$X_M=X_L\oplus C.\tag 5.1$$
If $L$ is simply connected, then there is a map $g: M\to L$ of
degree 1.
\endproclaim

\demo{\bf Proof} In the whole proof, the intersection forms are
understood in homology rings
defined by cap product.
We also use the same symbol for a surface and its homology class.
The proof is divided into three steps:

(1) First prove Proposition 5.1 in the case $M$ is also simply connected.

Suppose $X_M=X_L\oplus C$. Clearly $C$ is  symmetric and
unimodular, therefore $C$ represents the intersection form of a
simply connected 4-manifold $Q$ by Freedman's work [9]. So we have
$X_M=X_L\oplus X_Q=X_{L\#Q}$. By Whitehead's theorem [II. Theorem
2.1, 7], $M$ and $L\#Q$ are homotopy equivalent. Then there is a
homotopy equivalence $h : M\to L\# Q$, which is of degree one.
There is also an obvious degree one pinch $p: L\# Q\to L$. Then
$p\circ h: M\to L$ is a map of degree one.

(2) Then prove Proposition 5.1 in the case $X_M=X_L$.

Let $\Cal F=\{F_j, \, j=1,...,m\}$ be a set of oriented
embedded surfaces which provide
a basis of $\bar H_2(M)$.

Suppose $M$ is not simply connected (otherwise it is proved in
(1)). Let $\{ C_i$, $i=1,...,l\}$ be a set of disjoint simple
closed curves which generate $\pi_1(M)$ satisfying
$$(\cup_i N(C_i))  \cap (\cup_j F_j)=\emptyset,\tag 5.2$$

where $N(C_i)=C_i\times D^3$ is a regular neighborhood of $C_i$ in
$M$. Let $M^*$ be the simply connected 4-manifold obtained from
$M$ by surgery on $\{C_i\}$. Precisely

$$M^*=(M-\text{int} \cup _i C_i\times D^3)\cup _{\{h_i\}} (\cup _i D^2_i\times S^2),
\tag 5.3$$
where  the gluing map
$$h_i: S_i^1\times S^2=\partial D^2_i\times S^2\to
\partial (C_i \times D^3)=C_i\times S^2\tag 5.4$$
is an orientation reversing homeomorphism.

Since $\Cal F$ is linear independent in $\bar H_2(M)$, it is easy
to see that $\Cal F$ is still linear independent in $H_2(M^*)$.
Hence $\Cal F$ is a basis of a submodule $H$ of $H_2(M^*)$ and the
restriction of the intersection form $X_{M^*}$ on  $H$ is still
$X_M$, which   is isomorphic to $X_L$. In particular the
restriction of $X_{M^*}$ on $H$ is still unimodular. By [p.5, 8],
we have
$$X_{M^*}=X_M\oplus X_{H'}=X_L\oplus X_{H'},\tag 5.5$$

where $H'$ is the orthogonal complement of $H$ in $H_2(M^*)$, and

$$H_2(M^*)=H\oplus H'=H_2(L)\oplus H'.\tag 5.6$$
By part (1), we have a simply connected 4-manifold $Q$ with
$H_2(Q)=H'$ and a degree one map
$$\CD M^*@>h>> L\# Q@> p>> L,\endCD \tag 5.7$$
where $h$ is a homotopy equivalence and $p$ pinches $Q$. Let
$f=p\circ h$. For

$$ f_*: H_2(M^*)=\pi_2(M^*)\to H_2(L)=\pi_2(L),\tag 5.8$$

we have $H'=\text{Ker} f_* $.

Now we will transfer the  degree one map in (5.7) to a degree one map
$M\to L$.

Since $*_i\times S^2$ is disjoint from all $F_j$,  $*_i\times
S^2\in H'$ for each $i$, where $*_i\in D^2_i$. Then by (5.8)
$f(*_i\times S^2)$ is null homologous, and therefore is null
homotopic. Now we can  homotope $f$ so that $f(D^2_i\times S^2)=z$
for each $i$ by first shrinking $D^2_i\times S^2$ to its core
$*_i\times S^2$, then shrinking $*_i\times S^2$ to $z$.

Note $M^*-\text{int}\cup _i D^2_i\times S^2=M-\text{int} \cup _i
C_i\times D^3$. Since $f: M^*\to L$ is of degree one,
$f(D^2_i\times S^2)=z$ for all $i$, we have that the restriction

$$f| : (M-\text{int} \cup _i C_i\times D^3, \cup _i C_i\times S^2)\to
(L,z)\tag 5.9$$ is a map between pairs of degree 1.

Now we can extend $f|$ to $g: M\to L$ by sending each $C_i\times
D^3$ to $z$. Since $g$ sends $C_i\times D^3$ to $z$, the extension
$g:M\to L$ is of degree 1.

(3) Prove  Proposition 5.1 from (1) and (2).

By [9], there is a simply connected 4-manifold $Q$ with
$X_Q=X_M=X_L\oplus C$. Then there is a degree one map from $M$ to
$Q$ by (2) and a degree one map from $Q$ to $L$ by (1). Hence
there is a degree one map form $M$ to $L$. \qed

\demo{Proof of Theorem 3} The necessary part follows from Theorem
2, and the proof of the sufficient part is divided into two steps.

(1) Suppose first that $M$ is also simply connected. We may
suppose that $M$ and $L$ have presentations (3.1), and the bases
$\alpha$ and $\beta$ for $\bar H_2(M)$ and $\bar H_2(L)$ are
chosen as in (3.2). We may further suppose that $\alpha^*$ and
$\beta^*$ in Theorem 3 are the algebraic dual bases of $\alpha$
and $\beta$ respectively.

Then $f^*(\beta^*)=\alpha^* P$ implies that $f_*(\alpha)=\beta
P^\bot$, and by Lemma 2.4, the linking matrices for $M$ and $L$ in
the presentations (3.1) are $A$ and $B$ in Theorem 3 respectively.
Hence (1.4) in Theorem 1 holds. Since $n=2$, by Theorem 2.3 (2),
(1.5) is covered by (1.4). Hence there is a map $f: M\to L$ of
degree $k$ realizing $P$ by Theorem 1.

(2) If $M$ is not simply connected, let $Q$ be a simply connected
4-manifold with $X_Q=X_M$. By Proposition 5.1, there is a degree
one map $h: M\to Q$ which induces an isomorphism $h_*: \bar
H_2(M)\to H_2(Q)$, hence also induces an isomorphism $h^*:
H^2(Q)\to \bar H^2(M)$. Under the basis ${h^*}^{-1}(\alpha^*)$,
the intersection matrix of $X_Q$ is $A$. Now
${h^*}^{-1}(\alpha^*P)={h^*}^{-1}(\alpha^*)P$ and $P^\bot A P=kB$.

By step (1) there is a map $f': Q\to L$ of degree $k$ such that
${f'}^{*}\beta^*={h^*}^{-1}(\alpha^*P)$. And finally let
$f=f'\circ h$. $f$ is of degree $k$ and $f^*\beta^* =
\alpha^*P$.\qed

\vskip 0.5 true cm
{\bf \S 6. Proof of the corollaries.}

\demo{Proof of Corollary 1}
The first claim in Corollary 1 follows from Theorem 2 and the arguments
in the proof of Proposition 5.1 (1).

For the remaining,
since
$\beta_0(M)=\beta_0(L)=\beta_4(M)=\beta_4(L)=1$ and
$\beta_1(M)=\beta_1(L)=\beta_3(M)=\beta_3(L)=0$,
 $\chi(M)=\chi(L)$ implies that $\beta_2(M)=\beta_2(L)$, where
 $\beta_i$ is the $i$-th Betti number.
Since $f$ is of degree one, by Theorem 3, the intersection forms
of $M$ and $L$ must be isomorphic. By [II.Theorem 2.1, 7], $M$ and
$L$ are homotopy equivalent. (If the intersection form $X_M$ is
even, $M$ and $L$ are homeomorphic by [9]). \qed \vskip 0.5 true
cm

\demo{Proof of Corollary 2} It is known  that $\pi_{2n-1}(S^n)=
<\nu> \oplus G$, where $\nu$ has infinite order if $n$ is even and
is 0 if $n$ is odd, and $G$ is a finite group (see [p.318 and
p.325, 6]). Moreover if $n$ is even and $t\in \pi_{2n-1}(S^n)$,
then $t=\lambda H(t) \nu +\mu$, where $\mu \in G$, the Hopf
invariant $H(t)$ is the self-linking number defined by the map
$t:S^{2n-1}\to S^n$ and $\lambda=1$ if $n=2,4,8$ and $\lambda
=1/2$ otherwise (see [p.326, 6]).  It is also known that $H([s^n,
s^n])=2$ ([p.336, 6]).

Let $(A; t_1,...,t_m)$ be the complete homotopy invariant of $M$.

Now $t_v=\lambda H(t_v) \nu + \mu_v= a_{vv}\lambda\nu + \mu_v$ and
$[s^n, s^n]=2 \lambda\nu + \mu$, where $\mu_v, \mu\in G$. Let
$T(n)$ be the order of the torsion part of $\pi_{2n-1}(S^n)$ and
$k$ be a multiple of $2T(n)$ if $T(n)$ is even,  or a multiple of
$T(n)$ otherwise. Let $P=kI_m$, where $I_m$ is the $m$ by $m$ unit
matrix. Now we substitute all these information into each side of
(1.4) and (1.5) in Theorem 1. We first have
$$P^\bot A P=k^2 A.\tag 6.1$$

Now the right side of (1.5) in Theorem 1 can be reduced to

$$\sum_v p_{vr}t_v +(\sum _v \frac 12 p_{vr}(p_{vr}-1)a_{vv}+
\sum _{v<w}  p_{vr}p_{wr}a_{vw})[s^n,s^n]$$
$$=kt_r +\frac 12 k(k-1)a_{rr}[s^n,s^n]$$
$$=k(a_{rr}\lambda\nu+\mu_r) +\frac 12 k(k-1)a_{rr}(2\lambda\nu+\mu)$$
$$= k a_{rr}\lambda
\nu + k(k-1)a_{rr}\lambda \nu=k^2 a_{rr}\lambda\nu\tag 6.2$$ On
the other hand, clearly

$$k^2t_r=k^2( a_{rr}\lambda\nu+\mu_r)=k^2a_{rr}\lambda\nu. \tag 6.3$$

(6.1) implies that (1.4) of Theorem 1 holds for $k^2$, and (6.2)
and (6.3) imply that (1.5) of Theorem 1 holds for $k^2$. By
Theorem 1 there is a map of degree $k^2$ from $M$ to itself
realizing $P$.

In the case $n=2$, we have $T(2)=1$. We have finished the proof.
\qed

\vskip 0.5 true cm \demo{Proof of Corollary 3} Let $P$ be the
matrix realized by $f$. By (1.6) in Theorem 2 we have
$$P^\bot A P=kA,\tag 6.4$$
where both $A$, $P$ are square matrices of order $\beta_n(M)$. By
taking the determinants we have $|P|^2|A|=k^{\beta_n(M)}|A|$.
Thus $|P|^2=k^{\beta_n(M)}$. Since $\beta_n(M)$ is odd, $k$ itself
must be a perfect square. \qed

\vskip 0.5 true cm \demo{Proof of Corollary 4} From the equation
$P^\bot AP=kB$, i.e., (1.6) in Theorem 2, we have
$$kb_{ii}=\sum_{r,s} a_{rs}p_{ri}p_{si}
=\sum_{r} a_{rr}p_{ri}^2+ \sum_{r<s} (a_{rs}p_{ri}p_{si}+
a_{sr}p_{si}p_{ri}).\tag 6.5$$ Since $X_M$ is even, each $a_{rr}$
is even. Since $a_{rs}= \pm a_{sr}$, the right side of the
equation is even. Since $X_L$ is not even, $b_{ii}$ is odd for
some $i$. It follows $k$ is even.\qed \vskip 0.5 true cm

\demo{Proof of Corollary 5} If $M$ dominates $L$, then the rank of
the intersection form of $L$ is at most the rank of the
intersection form of $M$ by Theorem 2. For a given  positive
integer $m$, there are only finitely many non-isomorphic forms of
rank $m$ [p.18 and p.25, 8]. Since the torsion part of
$\pi_{2n-1}(S^n)$ is finite, the homotopy types of possible $L$
are finite. In the case of 4-manifolds, there are at most two
simply connected 4-manifolds with a given form $X$ [9]. So the
conclusion follows. \qed

\vskip 0.5 true cm

{\bf \S 7. Examples.}

Now we will apply the results above to some simple examples of 4-manifolds.

Let $F_g$, $CP^2$ and $T^4$ be the
surface of genus $g$, complex projective plane and 4-dimensional
torus respectively.

Let
$$
I_l=\oplus_l (1), \qquad
A_l=\oplus_l\pmatrix 0 & 1\\1 & 0\endpmatrix,\qquad
P_{l,k}=\oplus_l\pmatrix 0 & k\\1 & 0\endpmatrix,\qquad
B=\pmatrix 1 & 0\\0 & -1\endpmatrix.
$$

\vskip 0.3 true cm
{\bf Example 1.}
Any closed orientable 4-manifold $M$ with non-trivial intersection form
dominates $CP^2$.

{\bf Verification:} Suppose $M$ has intersection matrix
$A=(a_{ij})_{m\times m}$. If $a_{ii}\ne 0$ for some $i$, then let
$c_i=1$ and $c_r=0$ for $r\ne i$. If $a_{ii}=0$ for all $i$, then
pick any $a_{ij}\ne 0$, let $c_i=c_j=1$ and $c_r=0$ for $r\ne
i,j$. Let $C=(c_{i})_{m\times 1}$. Then
$$C^\bot A C=k(1).$$
Then $k=2a_{ij}$ if $i\ne j$
and $k=a_{ii}$ otherwise.
By Theorem 3, there is a map $f: M\to CP^2$ of degree $k\ne 0$.
\vskip 0.3 true cm

{\bf Example 2.} Suppose $M$ is a simply connected 4-manifold. (1)
If $P_{m\times 1}$ is realized  by a map $f: M\to CP^2$, then any
map $M\to CP^2$ realizing $P$ is homotopic to $f$. (2) If
$P_{m\times 2}$ is realized  by a map $f: M\to S^2\times S^2$,
then all maps $M\to S^2\times S^2$ realizing $P$ have 4 homotopy
classes. (3) There are infinitely many homotopy classes of maps
$(S^2\times S^2)\#(S^3\times S^1)\to S^2\times S^2$ realizing
$I_2$.

{\bf Verification:} By the exact sequence of fibration $S^5\to
CP^2$, we have $\pi_4(CP^2)=\pi_4(S^5)=0$. Moreover
$\pi_4(S^2\times S^2)= \pi_4(S^2)\oplus \pi_4(S^2)=Z_2\oplus Z_2$.
Then (1) and (2) follow from Theorem 1. For (3), consider the
composition

$$
(S^2\times S^2)\# (S^3\times S^1)
 \overset \, p_1\to \to
 (S^2\times S^2) V (S^3\times S^1)
\overset \, p_2\to \to
(S^2\times S^2) V S^3
\overset \, g_k\to \to
S^2\times S^2,
$$

where $p_1$ pinches the 3-sphere separating $S^2\times S^2$ and
$S^3\times S^1$ to a point to get the one point union of
$S^2\times S^2$ and $S^3\times S^1$, $p_2|S^3\times S^1 :
S^3\times S^1\to S^3$ is the projection and $p_2| S^2\times S^2$
is the identity, $g_k| S^3$ sends $S^3$ to the first factor of
$S^2\times S^2$ with Hopf invariant $k$ and $g_k|S^2\times S^2:
S^2\times S^2\to S^2\times S^2$ realizes $I_2$. Let $f_k= g_k\circ
p_2 \circ p_1$, then $\{f_k\}$ forms infinitely many homotopy
classes realizing $I_2$.

\vskip 0.3 true cm

{\bf Example 3.}
$D(CP^2\# (-CP^2), S^2\times S^2)=
D(S^2\times S^2, CP^2\#(-CP^2))=2Z$ and
$D(CP^2\# CP^2, S^2\times S^2)=
D(S^2\times S^2, CP^2\#CP^2)=
D(CP^2\# (-CP^2),  CP^2\#  CP^2)$ $=
D(CP^2\# CP^2, CP^2\# (-CP^2))=\{0\}$.

{\bf Verification:} Let $P= \pmatrix a & c\\b & d\endpmatrix$.
Since
$$P^\bot BP
=\pmatrix a^2-b^2 & ac-bd\\ac-bd & c^2-d^2\endpmatrix,$$ one can
verify easily that there is no solution $P$ for $P^\bot B P=k I_2$
for any non-zero integer $k$ and no solution $P$ for $P^\bot B P=k
A_1$  for any odd integer $k$. Moreover $P=\pmatrix 1 & k\\1 &
-k\endpmatrix$ is a solution for $P^\bot B P=2k A_1$.

Since $I_2$, $B$ and $A_1$ are the intersection forms of $CP^2\#
CP^2$, $CP^2\#(-CP^2)$ and $S^2\times S^2$ respectively, $D(CP^2\#
(-CP^2),  CP^2\#  CP^2)=\{0\}$ and $D(CP^2\# (-CP^2), S^2 $
$\times S^2)=2Z$ follow from Theorem 3.   The verification for the
remaining is  similar. \vskip 0.3 true cm

{\bf Example 4.} $D(T^4, \#_3 S^2\times S^2)=Z$.
That is for any $k$, there is a map of degree $k$ from 4-torus to the connected sum of 3 copies of $S^2\times S^2$.
In general $D(F_s\times F_r, \#_q S^2\times S^2)=Z$ if $q\le 2rs+1$.

{\bf Verification:}
 Note $T^4=S_1^1\times S_2^1\times S_3^1\times S_4^1$.
Let $a_i\in H^1(T^4)$ be the algebraic dual of $S_i^1$. Let
$c_{ij}=a_i\cup a_j$. Then under the basis $(c_{12}, c_{34},
c_{13}, c_{24}, c_{14}, c_{23})$ of $H^2(T^4)$, the intersection
matrix is $A_3$, which is the same as the matrix of $\#_3
S^2\times S^2$ under the obvious basis. Then it is easy to verify
$$P_{3,k}^\bot A_3P_{3,k} =k A_3.  $$
Hence there is a map $f_k: T^4\to \#_3 S^2\times S^2$ of degree $k$ by Theorem 3.

The verification of the second part is similar.
\vskip 0.3 true cm

{\bf REFERENCES.}

[1]  Edmonds, A.: {\it Deformation of maps to branched covering in
dimension 2,} Ann. Math. {\bf 110}:  113-125 (1979)

[2]  Wang,S.C.: Non-zero degree maps between 3-manifolds,
Proceedings of International Congress of Mathematicains, (Beijing
2002) Vol. II 457-468,  High Education Press, Beijing, 2002.

[3]  Duan H. B. and  Wang S.C.: {\it The degrees of maps between
manifolds,}  Math. Zeit. 244 67-89 (2003)

[4]  Dubrovin, B.A.,  Fomenko, A.T. and Novikov, S.: Modern
Geometry, II GTM Vol. 104 Springer-Verlag, 1984.

[5]  Wall, C.T.C.: {\it Classification of (n-1)-connected
2n-manifolds,} Ann. of Math. {\bf 75}: 163-198 (1962)

[6]  Hu, S.: Homotopy theory, Academic Press, New York and London,
1959.

[7]  Kirby, R.: The topology of 4-manifolds, Lecture Notes in Math
Vol. 1374, Springer-Verlag, 1989.

[8]  Milnor, J.  and  Husemoller, D.: Symmetric bilinear forms,
Springer-Verlag, Berlin and New York, 1973.

[9]  Freedman,  M.H.: {\it The topology of 4-dimensional
manifolds,} J. Diff. Geom. {\bf 17}: 357-453 (1982)

\bye